\documentclass{amsproc}
\usepackage{amsfonts,amssymb,amsmath,amsthm}
\usepackage{url}
\usepackage{enumerate}
\usepackage[all]{xypic}
\usepackage{mathrsfs}
\usepackage{graphicx}
\usepackage[matrix,arrow,curve,frame,arc]{xy}

\newtheorem{theorem}{Theorem}[section]

\newtheorem{corollary}[theorem]{Corollary}
\newtheorem{proposition}[theorem]{Proposition}

\theoremstyle{definition}

\newtheorem{example}[theorem]{Example}

\theoremstyle{remark}
\newtheorem{remark}[theorem]{Remark}

\newtheorem*{prob}{Problem}

\numberwithin{equation}{section}

\DeclareMathOperator{\mo}{mod}
\DeclareMathOperator{\ind}{ind}
\DeclareMathOperator{\rad}{rad}
\DeclareMathOperator{\Tr}{Tr}
\DeclareMathOperator{\pd}{pd}
\DeclareMathOperator{\id}{id}
\DeclareMathOperator{\gldim}{gl.dim}
\DeclareMathOperator{\Ext}{Ext}

\DeclareMathOperator{\Hom}{Hom}

\DeclareMathOperator{\End}{End}

\DeclareMathOperator{\add}{add}

\DeclareMathOperator{\Supp}{Supp}

\DeclareMathOperator{\ann}{ann}

\newcommand{\soc}{\mathrm{soc}}

\newcommand{\sC}{\mathscr{C}}

\newcommand{\sP}{\mathscr{P}}
\newcommand{\sQ}{\mathscr{Q}}
\newcommand{\sD}{\mathscr{D}}
\newcommand{\sT}{\mathscr{T}}

\newcommand{\cF}{\mathcal{F}}

\newcommand{\cA}{\mathcal{A}}
\newcommand{\Ab}{\cA b}
\newcommand{\bN}{\mathbb{N}}

\def\scr#1{{\scriptstyle{#1}}}

\author{Piotr Malicki and Andrzej Skowro\'nski}
\address{
Faculty of Mathematics and Computer Science\\
Nicolaus Copernicus University\\
Chopina 12/18, 87-100 Toru\'n, Poland}
\email{pmalicki@mat.umk.pl}
\email{skowron@mat.umk.pl}

\begin{document}

\title[Cycle-finite modules]{Cycle-finite modules over artin algebras}

\dedicatory{Dedicated to Jos\'e Antonio de la Pe\~na on the occasion of his 60th birthday}

\begin{abstract}
We describe the representation theory of finitely generated indecomposable modules over artin algebras which do not lie
on cycles of indecomposable modules involving homomorphisms from the infinite Jacobson radical of the module category.
\end{abstract}
\maketitle

\section{Introduction} \label{sect1}

Throughout the paper, by an algebra is meant an artin algebra over a fixed commutative artin ring $K$, which we shall assume (without loss of generality)
to be basic and indecomposable. For an algebra $A$, we denote by $\mo A$ the category of finitely generated right $A$-modules and by $\ind A$ the full
subcategory of $\mo A$ formed by the indecomposable modules. The Jacobson radical $\rad_A$ of $\mo A$ is the ideal generated by all nonisomorphisms between modules
in $\ind A$, and the infinite radical $\rad^{\infty}_A$ of $\mo A$ is the intersection of all powers $\rad^i_A$, $i\geq 1$, of $\rad_A$.
By a result of M. Auslander \cite{Au}, $\rad_A^{\infty}=0$ if and only if $A$ is of finite representation type, that is, $\ind A$ admits only a finite
number of pairwise nonisomorphic modules (see also \cite{KS0} for an alternative proof of this result).
On the other hand, if $A$ is of infinite representation type then $(\rad_A^{\infty})^2\neq 0$, by a result
proved in \cite{CMMS}.

An important combinatorial and homological invariant of the module category $\mo A$ of an algebra $A$ is its Auslander-Reiten quiver $\Gamma_A$.
Recall that $\Gamma_A$ is a valued translation quiver whose vertices are the isomorphism classes $\{X\}$ of modules $X$ in $\ind A$, the arrows correspond
to irreducible homomorphisms between modules in $\ind A$, and the translation is the Auslander-Reiten translation $\tau_A=D\Tr$.
We shall not distinguish between a module in $X$ in $\ind A$ and the corresponding vertex $\{X\}$ of $\Gamma_A$. If $A$ is an algebra of finite
representation type, then every nonzero nonisomorphism in $\ind A$ is a finite sum of composition of irreducible homomorphisms between modules in $\ind A$,
and hence we may recover $\mo A$ from the translation quiver $\Gamma_A$. In general, $\Gamma_A$ describes only the quotient category $\mo A/\rad^{\infty}_A$.

A prominent role in the representation theory of algebras is played by cycles of indecomposable modules (see \cite{MPS1}, \cite{MS4}, \cite{Ri1}, \cite{Sk7}).
Recall that a \textit{cycle} in $\ind A$ is a sequence
\[ M_0 \buildrel {f_1}\over {\hbox to 6mm{\rightarrowfill}} M_1 \to \cdots \to M_{r-1} \buildrel {f_r}\over {\hbox to 6mm{\rightarrowfill}} M_r=M_0 \]
of nonzero nonisomorphisms in $\ind A$ \cite{Ri1}, and such a cycle is said to be \textit{finite} if the homomorphisms $f_1,\ldots, f_r$
do not belong to $\rad_A^{\infty}$ (see \cite{AS1}, \cite{AS2}).
A module $M$ in $\ind A$ is said to be \textit{cycle-finite} if every cycle in $\ind A$ passing through $M$ is finite.
We note that this definition is more general then the one presented in \cite{MPS3}. Namely, every module $M$ in $\ind A$ which does not lie on a cycle
(\textit{directing module} in the sense of \cite{Ri1}) is cycle-finite.

If $A$ is an algebra of finite representation type, then all cycles in $\ind A$ are finite, and much of their properties is visible in the combinatorial
structure of the finite Auslander-Reiten quiver $\Gamma_A$ of $A$.
On the other hand, if $A$ is an indecomposable algebra of infinite representation type, then cycle-finite modules in $\ind A$ lie in infinite components
of the Auslander-Reiten quiver $\Gamma_A$ and usually belong to cycles containing arbitrary large number of pairwise nonisomorphic indecomposable modules.
Hence the study of cycle-finite modules over algebras of infinite representation type is a rather complicated problem.

The aim of this article is to present a rather complete representation theory of cycle-finite indecomposable modules over artin algebras.

For general results on the relevant representation theory we refer to the books \cite{ASS}, \cite{ARS}, \cite{Ri1}, \cite{SS1}, \cite{SS2}, \cite{SY3}, \cite{SY4}
and the survey articles \cite{CB2}, \cite{MPS1}, \cite{MS4}, \cite{Ri0}, \cite{Ri-86}, \cite{Sk7}.

\section{Preliminaries} \label{sect2}

Let $A$ be an algebra and $M$ a module in $\ind A$. An important information concerning the structure of $M$ is coded in the structure and properties of its
support algebra $\Supp(M)$ defined as follows. Consider a decomposition $A=P_M\oplus Q_M$ of $A$ in $\mo A$ such that the simple summands of the semisimple module
$P_M/\rad P_M$ are exactly the simple composition factors of $M$. Then $\Supp(M)=A/t_A(M)$, where $t_A(M)$ is the ideal in $A$ generated by the images of all
homomorphisms from $Q_{M}$ to $A$ in $\mo A$. We note that $M$ is an indecomposable module over $\Supp(M)$. Clearly, we may realistically hope to describe the
structure of $\Supp(M)$ only for modules $M$ having some distinguished properties.
For example, if $M$ is a directing module in $\ind A$, then the support algebra $\Supp(M)$
of $M$ over an algebra $A$ is a tilted algebra $\End_H(T)$, for a hereditary algebra $H$ and a tilting module $T$ in $\mo H$, and $M$ is
isomorphic to the image $\Hom_H(T,I)$ of an indecomposable injective module $I$ in $\mo H$ via the functor $\Hom_H(T,-): \mo H \to \mo\End_H(T)$
(see \cite{Ri1} and \cite{JMS1}, \cite{JMS2} for the corresponding result over arbitrary artin algebra).

Let $A$ be an algebra and $\Gamma_A$ be the Auslander-Reiten quiver of $A$.
By a component of $\Gamma_A$ we mean a connected component of the quiver $\Gamma_A$.
For a~component $\sC$ of $\Gamma_A$, we denote by $\ann_A(\sC)$ the annihilator of $\sC$ in $A$, that is, the intersection of the annihilators
$\{a\in A\mid Ma=0\}$ of all modules $M$ in $\sC$, and by $B(\sC)$ the quotient algebra $A/\ann_A(\sC)$, called the \textit{faithful algebra} of $\sC$.
We note that $\sC$ is a~faithful component of $\Gamma_{B(\sC)}$.
For a module $M$ in $\Gamma_A$, the $\tau_A$-\textit{orbit} of $M$ is the set $\mathcal{O}(M)$ of all possible vertices in $\Gamma_A$ of the form
$\tau_A^nM$, $n\in\mathbb{Z}$. The $\tau_A$-orbit $\mathcal{O}(M)$ of $M$ (respectively, the module $M$) is said to be \textit{periodic} if
$M\cong\tau_A^nM$ for some $n\geq 1$. A module $M$ in $\Gamma_A$ is said to be an \textit{acyclic module} if $M$ does not lie on an oriented cycle in
$\Gamma_A$, and otherwise a \textit{cyclic module}.
A component $\sC$ of $\Gamma_A$ without oriented cycles is said to be \textit{acyclic}.
Dually, a component $\sC$ of $\Gamma_A$ is said to be \textit{cyclic} if any module in $\sC$ lies on an oriented cycle of $\sC$.
Following \cite{MS1}, we denote by $_c\sC$ the full translation subquiver of $\sC$ obtained by removing all acyclic modules and the arrows attached to them,
and call it the \textit{cyclic part} of $\sC$. The connected translation subquivers of $_c\sC$ are said to be \textit{cyclic components} of $\sC$.
It was shown in \cite[Proposition 5.1]{MS1} that two modules $M$ and $N$ in $_c\sC$ belong to the same cyclic component of $\sC$ if there is an~oriented cycle in $\sC$
passing through $M$ and $N$.
All the modules in a cyclic component $\sC$ of $\Gamma_A$ containing a~cyclic-finite module are cyclic-finite, and such an cyclic
component $\sC$ is said to be a~\textit{cycle-finite cyclic component}.
For a subquiver $\sC$ of $\Gamma_A$, we consider a~decomposition $A=P_{\sC}\oplus Q_{\sC}$ of $A$ in $\mo A$ such that the simple summands
of the semisimple module $P_{\sC}/\rad P_{\sC}$ are exactly the simple composition factors of indecomposable modules in $\sC$, the ideal $t_A(\sC)$
in $A$ generated by the images of all homomorphisms from $Q_{\sC}$ to $A$ in $\mo A$, and call the quotient algebra $\Supp(\sC)=A/t_A(\sC)$ the
\textit{support algebra} of $\sC$.
Let $M$ be a nondirecting cycle-finite module in $\ind A$.
Observe that $M$ belongs to a unique cyclic component $\sC(M)$ of $\Gamma_A$ consisting entirely of cycle-finite
indecomposable modules, and the support algebra $\Supp(M)$ of $M$ is a quotient algebra of the support algebra $\Supp(\sC(M))$ of $\sC(M)$.
Moreover, by a result stated in \cite[Corollary 1.3]{MPS3}, the support algebra $\Supp(\sC)$ of a cycle-finite cyclic component $\sC$ of $\Gamma_A$ is isomorphic
to an algebra of the form $e_{\sC}Ae_{\sC}$ for an idempotent $e_{\sC}$ of $A$ whose primitive summands correspond to the vertices of a convex subquiver of the valued quiver
$Q_A$ of $A$. On the other hand, the support algebra $\Supp(M)$ of a cycle-finite module $M$ in $\ind A$ is not necessarily an algebra of the form
$eAe$ for an idempotent $e$ of $A$ (see \cite[Section 6]{MPS3}).
%
%
A component $\sC$ of $\Gamma_A$ is called \textit{regular} if $\sC$ contains neither a projective module nor an injective module,
and \textit{semiregular} if $\sC$ does not contain both a projective and an injective module. It has been shown in \cite{Li1}
and \cite{Zh} that a regular component $\sC$ of $\Gamma_A$ contains an oriented cycle if and only if $\sC$ is a \textit{stable tube}
(is of the form ${\Bbb Z}{\Bbb A}_{\infty}/(\tau^{r})$, for a positive integer $r$). Moreover, S. Liu proved in \cite{Li2} that a semiregular component
$\sC$ of $\Gamma_A$ contains an oriented cycle if and only if $\sC$ is a \textit{ray tube} (obtained from a stable tube by a finite number
(possibly zero) of ray insertions) or a \textit{coray tube} (obtained from a stable tube by a finite number (possibly zero) of coray insertions).
A component $\sC$ of $\Gamma_A$ is called \textit{postprojective} if $\sC$ is acyclic and every module in $\sC$
belongs to the $\tau_A$-orbit of a projective module. Dually, a component $\sC$ of $\Gamma_A$ is called \textit{preinjective} if $\sC$ is acyclic
and every module in $\sC$ belongs to the $\tau_A$-orbit of an injective module.
An indecomposable module $X$ in a component $\sC$ of $\Gamma_A$ is said to be \emph{right coherent} if there is in $\sC$ an infinite sectional path
\[X = X_1 \longrightarrow X_2 \longrightarrow \cdots \longrightarrow X_i\longrightarrow X_{i+1} \longrightarrow X_{i+2} \longrightarrow \cdots \]
Dually, an indecomposable module $Y$ in $\sC$ is said to be \emph{left coherent} if there is in $\sC$ an infinite sectional path
\[ \cdots \longrightarrow Y_{j+2} \longrightarrow Y_{j+1} \longrightarrow Y_j \longrightarrow \cdots \longrightarrow Y_2 \longrightarrow Y_1 = Y. \]
A~module $Z$ in $\sC$ is said to be \emph{coherent} if $Z$ is left and right coherent.
A component $\sC$ of $\Gamma_A$ is said to be \textit{coherent} \cite{MS1} (see also \cite{DR}) if every projective module in $\sC$ is right coherent
and every injective module in $\sC$ is left coherent.
Further, a component $\sC$ of $\Gamma_A$ is said to be \textit{almost cyclic} if its cyclic part $_c{\sC}$ is a~cofinite subquiver of $\sC$.
We note that the stable tubes, ray tubes and coray tubes of $\Gamma_A$ are special types of almost cyclic coherent components.
In general, it has been proved in \cite{MS1} that a component $\sC$ of $\Gamma_A$ is almost cyclic and coherent
if and only if $\sC$ is a \textit{generalized multicoil}, obtained from a finite family of stable tubes by a sequence of admissible operations (ad~1)-(ad~5)
and their duals (ad~1$^*$)-(ad~5$^*$).
We refer to \cite[Section 2]{MS1} for a~detailed description of these admissible operations and generalized multicoils.
In particular, one knows that all arrows of a~generalized multicoil have trivial valuation.
On the other hand, a component $\sC$ of $\Gamma_A$ is said to be \textit{almost acyclic} if all but finitely many
modules of $\sC$ are acyclic. It has been proved by I. Reiten and the second named author in \cite{RS2} that a component $\sC$ of $\Gamma_A$ is almost acyclic
if and only if $\sC$ admits a multisection $\Delta$. Moreover, for an almost acyclic component $\sC$ of $\Gamma_A$, there exists a finite convex
subquiver $c(\sC)$ of $\sC$ (possibly empty), called the \textit{core} of $\sC$, containing all modules lying on oriented cycles in $\sC$
(see \cite{RS2} for details).

Let $A$ be an algebra.
A family ${\sC}$ = $({\sC}_{i})_{i \in I}$ of components of $\Gamma_A$ is said to be \textit{generalized standard} if $\rad_A^{\infty}(X,Y)=0$ for all modules
$X$ and $Y$ in $\sC$ \cite{Sk3}, and \textit{sincere} if every simple module in $\mo A$ occurs as a composition factor of a module in $\sC$.
Two components $\sC$ and $\sD$ of an Auslander-Reiten quiver $\Gamma_A$ are said to be \textit{orthogonal} if $\Hom_A(X,Y) = 0$ and $\Hom_A(Y,X) = 0$
for all modules $X \in \sC$ and $Y \in \sD$. We also note that if $\sC$ and $\sD$ are distinct components of $\Gamma_A$ then $\Hom_A(X,Y) = \rad_A^{\infty}(X,Y)$
for all modules $X \in \sC$ and $Y \in \sD$.
Observe that a family ${\sC}$ = $({\sC}_{i})_{i \in I}$ of components of $\Gamma_A$ is generalized standard if and only if the components
$\sC_i$, $i\in I$, are generalized standard and pairwise orthogonal.
A prominent role in the representation theory of algebras is played by the algebras with separating families of Auslander-Reiten components.
A~concept of a~separating family of tubes has been introduced by C. M. Ringel in \cite{Ri0}, \cite{Ri1} who proved that they occur in the Auslander-Reiten quivers
of hereditary algebras of Euclidean type, tubular algebras, and canonical algebras.
More generally, following I. Assem, A. Skowro\'nski and B. Tom\'e \cite{AST}, a family ${\sC}$ = $({\sC}_{i})_{i \in I}$ in $\Gamma_A$ is said to be
\textit{separating} in $\mo A$ if the components of $\Gamma_A$ split into three disjoint families ${\sP}^A$, ${\sC}^A={\sC}$ and ${\sQ}^A$ such that
the following conditions are satisfied:
\begin{enumerate}
\item[(S1)] ${\sC}^A$ is a sincere generalized standard family of components;
\item[(S2)] $\Hom_{A}({\sQ}^A,{\sP}^A) = 0$, $\Hom_{A}({\sQ}^A,{\sC}^A)=0$, $\Hom_{A}({\sC}^A,{\sP}^A) = 0$;
\item[(S3)] any homomorphism from ${\sP}^A$ to ${\sQ}^A$ in $\mo A$ factors through the additive category $\add({\sC}^A)$ of ${\sC}^A$.
\end{enumerate}
\noindent Then we say that ${\sC}^A$ separates ${\sP}^A$ from ${\sQ}^A$ and write
\[\Gamma_A={\sP}^A \cup {\sC}^A \cup {\sQ}^A.\]
We note that then ${\sP}^A$ and ${\sQ}^A$ are uniquely determined by ${\sC}^A$ (see \cite[(2.1)]{AST} or \cite[(3.1)]{Ri1}).
Moreover, we have $\ann_A({\sC}^A)=0$, so ${\sC}^A$ is a~faithful family of components of $\Gamma_A$.
We note that if $A$ is an algebra of finite representation type that ${\sC}^A = \Gamma_A$ is trivially a~unique separating component of $\Gamma_A$,
with ${\sP}^A$ and ${\sQ}^A$ being empty.

\section{Quasitilted algebras} \label{sect-qua}

In the representation theory of algebras an important role is played by the canonical algebras introduced by C. M. Ringel in \cite{Ri1} and \cite{Ri2}.
We refer to \cite[Appendix]{Ri2} for an elementary definition of the canonical algebras proposed by W. Crawley-Boevey.
Every canonical algebra $\Lambda$ is of global dimension at most two. Moreover, the following theorem due to C. M. Ringel \cite{Ri2} describes the general shape
of the Auslander-Reiten quiver of a canonical algebra.
\begin{theorem}
\label{t65}%
Let $\Lambda$ be a canonical algebra. Then the general shape
of the Auslander-Reiten quiver $\Gamma_{\Lambda}$ of $\Lambda$ is as follows
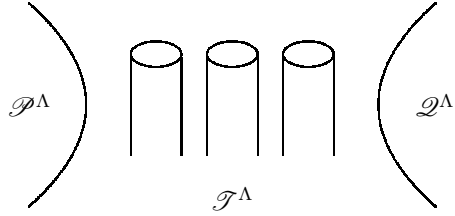
\[
\xymatrix@M=0pc@L=0pc@C=.8pc@R=1.6pc{
  \ar@/^5ex/@{-}[dddd] &&&&  && & && & &&  &&&& \ar@/_5ex/@{-}[dddd]  \\
  &&&&
  \ar@{-}[dd] & {}\save[] *{\xycircle<.8pc,.4pc>{}} \restore & \ar@{-}[dd] &
  \ar@{-}[dd] & {}\save[] *{\xycircle<.8pc,.4pc>{}} \restore & \ar@{-}[dd] &
  \ar@{-}[dd] & {}\save[] *{\xycircle<.8pc,.4pc>{}} \restore & \ar@{-}[dd] \\
  {}\save[] *{\sP^{\Lambda}} \restore &&&&  && & && & &&  &&&&
  {}\save[] *{\sQ^{\Lambda}} \restore \\
  &&&& && & && & &&&& \\
  &&&& && & & {}\save[] +<0mm,1mm> *{\sT^{\Lambda}} \restore & && & &&&& \\
}
\]
where $\sP^{\Lambda}$ is a family of components containing a unique postprojective component $\sP(\Lambda)$ and all indecomposable projective $\Lambda$-modules,
$\sQ^{\Lambda}$ is a family of components containing a unique preinjective component $\sQ(\Lambda)$ and all indecomposable injective $\Lambda$-modules,
and $\sT^{\Lambda}$ is an infinite family of faithful generalized standard stable tubes separating $\sP^{\Lambda}$ from $\sQ^{\Lambda}$.
In particular, we have $\pd_{\Lambda} X \leq 1$ for all modules $X$ in $\sP^{\Lambda} \cup \sT^{\Lambda}$
and $\id_{\Lambda} Y \leq 1$ for all modules $Y$ in $\sT^{\Lambda} \cup \sQ^{\Lambda}$.
\end{theorem}
Let $\Lambda$ be a canonical algebra. An algebra $C$ of the form $\End_{\Lambda}(T)$, where $T$ is a tilting $\Lambda$-module from the additive category
$\add(\sP^{\Lambda})$ of $\sP^{\Lambda}$ is called a \textit{concealed canonical algebra of type $\Lambda$}.
Then we have the following theorem describing the general shape of the Auslander-Reiten quiver of a concealed canonical algebra, which is
a consequence of Theorem \ref{t65} and the tilting theory.

\begin{theorem}
\label{t66}%
Let $\Lambda$ be a canonical algebra, $T$ a tilting module in $\add(\sP^{\Lambda})$, and $C=\End_{\Lambda}(T)$
the associated concealed canonical algebra. Then the general shape of the Auslander-Reiten quiver $\Gamma_C$ of $C$ is as follows
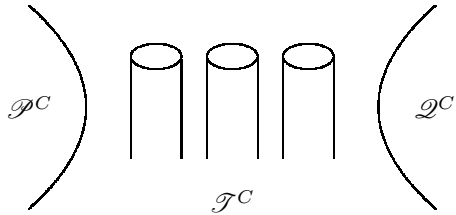
\[
\xymatrix@M=0pc@L=0pc@C=.8pc@R=1.6pc{
  \ar@/^5ex/@{-}[dddd] &&&&  && & && & &&  &&&& \ar@/_5ex/@{-}[dddd]  \\
  &&&&
  \ar@{-}[dd] & {}\save[] *{\xycircle<.8pc,.4pc>{}} \restore & \ar@{-}[dd] &
  \ar@{-}[dd] & {}\save[] *{\xycircle<.8pc,.4pc>{}} \restore & \ar@{-}[dd] &
  \ar@{-}[dd] & {}\save[] *{\xycircle<.8pc,.4pc>{}} \restore & \ar@{-}[dd] \\
  {}\save[] *{\sP^{C}} \restore &&&&  && & && & &&  &&&&
  {}\save[] *{\sQ^{C}} \restore \\
  &&&& && & && & &&&& \\
  &&&& && & & {}\save[] +<0mm,1mm> *{\sT^{C}} \restore & && & &&&& \\
}
\]
where $\sP^{C}$ is a family of components containing a unique postprojective component $\sP(C)$ and all indecomposable projective $C$-modules,
$\sQ^{C}$ is a family of components containing a unique preinjective component $\sQ(C)$ and all indecomposable injective $C$-modules,
$\sT^{C} = \Hom_{\Lambda}(T, \sT^{\Lambda})$ is an infinite family of faithful pairwise orthogonal generalized standard stable tubes separating $\sP^{C}$ from $\sQ^{C}$.
In particular, we have $\pd_{C} X \leq 1$ for all modules $X$ in $\sP^{C} \cup \sT^{C}$,
$\id_{C} Y \leq 1$ for all modules $Y$ in $\sT^{C} \cup \sQ^{C}$, and $\gldim C \leq 2$.
\end{theorem}
The following characterization of concealed canonical algebras has been established by J. A. de la Pe\~na and H. Lenzing in \cite{LP2}.

\begin{theorem}
Let $A$ be an algebra. The following statements are equivalent.
\begin{enumerate}
\item[\rm(1)] $A$ is a concealed canonical algebra.
\item[\rm(2)] $\Gamma_A$ admits a separating family $\sT^A$ of stable tubes.
\end{enumerate}
\end{theorem}
The concealed canonical algebras form a distinguished class of \textit{quasitilted algebras}, which are the endomorphism algebras $\End_{\mathscr H}(T)$ of tilting
objects $T$ in abelian hereditary $K$-categories ${\mathscr H}$ \cite{HRS}.
The following characterization of quasitilted algebras has been established by D. Happel, I. Reiten and S. O. Smal{\o} in \cite{HRS}.
\begin{theorem}
\label{th:3.7}%
Let $A$ be an algebra. The following statements are equivalent.
\begin{enumerate}
\item[\rm(1)] $A$ is a quasitilted algebra.
\item[\rm(2)] $\gldim A \leq 2$ and every module $X$ in $\ind A$ satisfies $\pd_AX\leq 1$ or $\id_AX\leq 1$.
\end{enumerate}
\end{theorem}
Recall briefly that the \textit{tilted algebras} are the endomorphism algebras $\End_{H}(T)$ of tilting modules $T$ over hereditary algebras $H$
and the \textit{quasitilted algebras of canonical type} are the endomorphism algebras $\End_{\mathscr H}(T)$ of tilting objects $T$ in abelian hereditary categories ${\mathscr H}$
whose derived category $D^b({\mathscr H})$ is equivalent to the derived category $D^b(\mo\Lambda)$ of the module category $\mo\Lambda$ of a canonical algebra $\Lambda$.
The next result proved in \cite{HRe} shows that there are only two
classes of quasitilted algebras.

\begin{theorem}
\label{t616}%
Let $A$ be a quasitilted algebra. Then $A$ is either a tilted algebra or a quasitilted algebra of canonical type.
\end{theorem}
The structure of representation-infinite quasitilted algebras of canonical type has been described by H. Lenzing and A. Skowro\'nski in \cite{LS1}
(see also \cite{S7}). In \cite{LS1} the concept of a semiregular branch enlargement of a concealed canonical algebra has been introduced and
the following characterization of quasitilted algebras of canonical type established,
which extends characterizations of concealed canonical algebras proved in \cite{LP2}, \cite{RS} and \cite{S6}.

\begin{theorem}
Let $A$ be an algebra. The following statements are equivalent.
\begin{enumerate}
  \item[\rm(1)] $A$ is a representation-infinite quasitilted algebra of canonical type.
  \item[\rm(2)] $A$ is a semiregular branch enlargement of a concealed canonical algebra.
  \item[\rm(3)] $\Gamma_A$ admits a separating family $\sT^A$ of semiregular tubes.
\end{enumerate}
\end{theorem}

Here, by a semiregular tube we mean a ray tube or a coray tube.

We may visualise the shape of the Auslander-Reiten quiver $\Gamma_A$ of
a representation-infinite quasitilted algebra $A$ of canonical type as follows
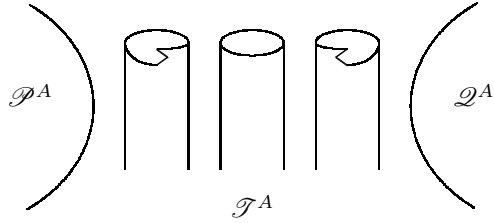
\[
\xy
0;/r1pc/:0
,{\ellipse(1,.4),=:a(180){-}}
,{\ellipse(1,.2)_,=:a(90){-}}
,{\ellipse(1,.7)__,=:a(90){-}}
*\dir{}="a",;p+(0,-.2)*\dir{}="b",**\dir{},"b"
;p+(.35,-.25)*\dir{}="b",**\dir{-},"b"
;p+(-.35,-.25)*\dir{}="b",**\dir{-},"b"
;p+(-1,.7)
;p+(0,-4)*\dir{}="b",**\dir{-},"b"
;p+(2,0)
;p+(0,4)*\dir{}="b",**\dir{-},"b"
;p+(2,0)
,{\ellipse(1,.4){-}}
;p+(-1,0)
;p+(0,-4)*\dir{}="b",**\dir{-},"b"
;p+(2,0)
;p+(0,4)*\dir{}="b",**\dir{-},"b"
;p+(2,0)
,{\ellipse(1,.4),=:a(180){-}}
,{\ellipse(1,.2)__,=:a(90){-}}
,{\ellipse(1,.7)_,=:a(90){-}}
*\dir{}="a",;p+(0,-.2)*\dir{}="b",**\dir{},"b"
;p+(-.35,-.25)*\dir{}="b",**\dir{-},"b"
;p+(.35,-.25)*\dir{}="b",**\dir{-},"b"
;p+(-1,.7)
;p+(0,-4)*\dir{}="b",**\dir{-},"b"
;p+(2,0)
;p+(0,4)*\dir{}="b",**\dir{-},"b"
;p+(6,-2)
,{\ellipse(5,4):a(120),,=:a(120){-}}
*\dir{}="a",;p+(-20,0)*\dir{}="b",**\dir{},"b"
,{\ellipse(5,4):a(210),,=:a(120){-}}
;p+(3,0)*+!{\sP^A}
;p+(14,0)*+!{\sQ^A}
;p+(-7,-3.5)*+!{\sT^A}
\endxy
\]
where $\sT^A$ is a faithful infinite family of pairwise orthogonal, generalized standard ray or coray tubes, separating $\sP^A$ from $\sQ^A$.
Then every indecomposable projective $A$-module lies in $\sP^A \cup \sT^A$ while every indecomposable injective $A$-module lies in $\sT^A \cup \sQ^A$.

The following example from \cite{MS4} illustrates the above theorem.
\begin{example}
Let $K$ be an algebraically closed field, $Q$ the quiver
$$
  
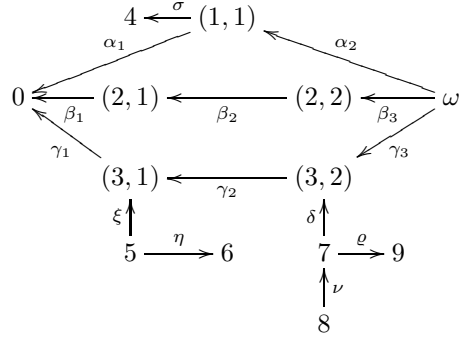
\xymatrix@C=0.7pc@R=1.2pc{
    && 4 & (1,1) \ar[l]_{\sigma} \ar[llld]_{\alpha_1} \\
    0 && (2,1) \ar[ll]^{\beta_1} && (2,2) \ar[ll]^{\beta_2} &&
         \omega \ar[lld]^{\gamma_3} \ar[ll]^{\beta_3} \ar[lllu]_{\alpha_2} \\
    && (3,1) \ar[llu]^{\gamma_1} && (3,2) \ar[ll]^{\gamma_2} \\
    && 5 \ar[u]^{\xi} \ar[r]^{\eta} & 6 &
      7 \ar[u]^{\delta} \ar[r]^{\varrho} & 9 \\
    &&&& 8 \ar[u]_{\nu}
  }
$$
$I$ the ideal of $K Q$ generated by the elements
$\alpha_2 \alpha_1 + \beta_3 \beta_2 \beta_1 + \gamma_3 \gamma_2 \gamma_1$, $\alpha_2\sigma$, $\xi \gamma_1$, $\delta \gamma_2$, $\nu \varrho$,
and $A=KQ/I$ the associated bound quiver algebra. Then $A$ is a representation-infinite quasitilted algebra of canonical type.
Moreover, $A$ is a semiregular branch enlargement of the canonical algebra $C = K\Delta / J$, where $\Delta$ is the full subquiver of $Q$ given
by the vertices $0$, $\omega$, $(1,1)$, $(2,1)$, $(2,2)$, $(3,1)$, $(3,2)$, and $J$ is the ideal of $K \Delta$ generated by
$\alpha_2 \alpha_1 + \beta_3 \beta_2 \beta_1 + \gamma_3 \gamma_2 \gamma_1$. The Auslander-Reiten quiver
$\Gamma_A$ of $A$ has a disjoint union form
\[ \Gamma_A = \sP^A \cup \sT^A \cup \sQ^A \]
where $\sT^A$ a family of semiregular tubes separating $\sP^A$ from $\sQ^A$. Moreover,
$\sT^A$ is a ${\Bbb P}_1(K)$-family $(\sT^A_{\lambda})_{\lambda\in {\Bbb P}_1(K)}$ of semiregular tubes consisting of:
\begin{itemize}
\item a coray tube $\sT_{\infty}^A$ obtained from the stable tube $\sT_{\infty}^C$ of $\Gamma_C$ of rank $2$ by one coray insertion,
\item a stable tube $\sT_0^A = \sT_0^C$ of rank $3$ formed by indecomposable $C$-modules,
\item a ray tube $\sT_1^A$ obtained from the stable tube $\sT_{1}^C$ of $\Gamma_C$ of rank $3$ by $5$ ray insertions,
\item the infinite family $\sT^A_{\lambda}, \lambda\in {\Bbb P}_1(K)\setminus\{\infty,0,1\}$, of stable tubes of rank $1$,
consisting of indecomposable $C$-modules.
\end{itemize}
We refer to \cite[Example 6.19]{MS4} for a detailed description of the families $\sP^A, \sT^A$ and $\sQ^A$.
\end{example}

\section{Cycle-finite cyclic Auslander-Reiten components} \label{sect4}

In this section we describe the support algebras of cycle-finite cyclic components of the Auslander-Reiten quivers of artin algebras.
The description splits into two cases.
In the case when a cycle-finite cyclic component $\sC$ of $\Gamma_A$ is infinite, the support algebra $\Supp(\sC)$ is a suitable gluing of finitely many
generalized multicoil algebras (introduced by the authors in \cite{MS2}) and algebras of finite representation type, and $\sC$ is the
corresponding gluing of the associated cyclic generalized multicoils via finite translation quivers (Theorem \ref{thm-1-1}).
In the second case when a cycle-finite cyclic component $\sC$ is finite, the support algebra $\Supp(\sC)$ is a generalized double tilted algebra
(in the sense of I. Reiten and A. Skowro\'nski \cite{RS2}) and $\sC$ is the core of the connecting component of this algebra (Theorem \ref{thm-1-2}).

In order to present the first case, we need the class of generalized multicoil algebras.
Recall that following \cite{MS2}, an algebra $A$ is called a \textit{generalized multicoil algebra}, if $A$ is a generalized multicoil enlargement of a product
$C=C_1\times\ldots\times C_m$ of concealed canonical algebras $C_1, \ldots, C_m$ using modules from the separating families
${\sT}^{C_1}, \ldots, {\sT}^{C_m}$ of stable tubes of $\Gamma_{C_1}, \ldots, \Gamma_{C_m}$ and a sequence of admissible operations of types
(ad~1)-(ad~5) and their duals (ad~1$^*$)-(ad~5$^*$). The following result has been established in \cite[Theorem A]{MS2}.
\begin{theorem} \label{thm-41}
Let $A$ be an algebra. The following statements are equivalent.
\begin{enumerate}
\item[\rm(1)] $A$ is a generalized multicoil algebra.
\item[\rm(2)] $\Gamma_A$ admits a separating family of almost cyclic coherent components.
\end{enumerate}
\end{theorem}

The following consequence of \cite[Theorem C]{MS2} describes the structure of the Auslander-Reiten quivers of generalized multicoil algebras.

\begin{theorem} \label{thm-43}
Let $A$ be a generalized multicoil algebra obtained from a family $C_1, \ldots, C_m$ of concealed canonical algebras. Then there are unique quotient algebras
$A^{(l)}$ and $A^{(r)}$ of $A$ such that the following statements hold:
\begin{enumerate}
\item[\rm(1)] $A^{(l)}$ is a product of quasitilted algebras of canonical type having separating families of coray tubes.
\item[\rm(2)] $A^{(r)}$ is a product of quasitilted algebras of canonical type having separating families of ray tubes.
\item[\rm(3)] The Auslander-Reiten quiver $\Gamma_A$ has a disjoint union decomposition
\[ \Gamma_A = \sP^A \cup \sC^A \cup \sQ^A, \]
where
\begin{enumerate}
\item[\rm(a)] ${\sP}^A$ is the left part ${\sP}^{A^{(l)}}$ in a decomposition $\Gamma_{A^{(l)}}={\sP}^{A^{(l)}}\cup {\sT}^{A^{(l)}}\cup {\sQ}^{A^{(l)}}$ of the Auslander-Reiten quiver
$\Gamma_{A^{(l)}}$ of the algebra $A^{(l)}$, with ${\sT}^{A^{(l)}}$ a family of coray tubes separating ${\sP}^{A^{(l)}}$ from ${\sQ}^{A^{(l)}}$;
\item[\rm(b)] ${\sQ}^A$ is the right part ${\sQ}^{A^{(r)}}$ in a decomposition $\Gamma_{A^{(r)}}={\sP}^{A^{(r)}}\cup {\sT}^{A^{(r)}}\cup {\sQ}^{A^{(r)}}$ of the Auslander-Reiten quiver
$\Gamma_{A^{(r)}}$ of the algebra $A^{(r)}$, with ${\sT}^{A^{(r)}}$ a family of ray tubes separating ${\sP}^{A^{(r)}}$ from ${\sQ}^{A^{(r)}}$;
\item[\rm(c)] ${\sC}^A$ is a family of generalized multicoils separating ${\sP}^A$ from ${\sQ}^A$, obtained from stable tubes in the
separating families ${\sT}^{C_1}, \ldots, {\sT}^{C_m}$ of stable tubes of the Auslander-Reiten quivers $\Gamma_{C_1}, \ldots, \Gamma_{C_m}$
of the concealed canonical algebras $C_1, \ldots, C_m$ by a sequence of admissible operations of types (ad~1)-(ad~5) and their duals (ad~1$^*$)-(ad~5$^*$),
corresponding to the admissible operations leading from $C=C_1\times\ldots\times C_m$ to $A$;
\item[\rm(d)] ${\sC}^A$ consists of cycle-finite modules and contains all indecomposable modules of ${\sT}^{A^{(l)}}$ and ${\sT}^{A^{(r)}}$;
\item[\rm(e)] ${\sP}^A$ contains all indecomposable modules of ${\sP}^{A^{(r)}}$;
\item[\rm(f)] ${\sQ}^A$ contains all indecomposable modules of ${\sQ}^{A^{(l)}}$.
\end{enumerate}
\end{enumerate}
\end{theorem}
Moreover, in the above notation, we have the following consequences of \cite[Theorem E]{MS2} describing the homological properties of modules over generalized multicoil algebras
\begin{itemize}
\item $\gldim A\leq 3$;
\item $\pd_AX\leq 1$ for any indecomposable module $X$ in ${\sP}^A$;
\item $\id_AY\leq 1$ for any indecomposable module $Y$ in ${\sQ}^A$;
\item $\pd_AM\leq 2$ and $\id_AM\leq 2$ for any indecomposable module $M$ in ${\sC}^A$.
\end{itemize}

The algebra $A^{(l)}$ is said to be the \textit{left quasitilted algebra} of $A$ and the algebra $A^{(r)}$ is said to be the \textit{right quasitilted algebra} of $A$.

The following theorem from \cite[Theorem 1.1]{MPS3} describes the support algebras of infinite cycle-finite cyclic components of the Auslander-Reiten quivers
of artin algebras.
\begin{theorem} \label{thm-1-1}
Let $A$ be an algebra and $\sC$ be an infinite cycle-finite component of $_c\Gamma_A$. Then there exist infinite
full translation subquivers $\sC_1, \ldots, \sC_r$ of $\sC$ such that the following statements hold.
\begin{enumerate}
\item[\rm(1)] For each $i\in\{1,\ldots,r\}$, $\sC_i$ is a cyclic coherent full translation subquiver of $\Gamma_A$.
\item[\rm(2)] For each $i\in\{1,\ldots,r\}$, $\Supp(\sC_i)=B(\sC_i)$ and is a generalized multicoil algebra.
\item[\rm(3)] $\sC_1, \ldots, \sC_r$ are pairwise disjoint full translation subquivers of $\sC$ and $\sC^{cc}=\sC_1\cup\ldots\cup\sC_r$ is a maximal cyclic coherent
and cofinite full translation subquiver of $\sC$.
\item[\rm(4)] $B(\sC\setminus\sC^{cc})$ is of finite representation type.
\item[\rm(5)] $\Supp(\sC)=B(\sC)$.
\end{enumerate}
\end{theorem}
\begin{remark}
It follows from the above theorem that all but finitely many modules lying in an infinite cycle-finite component $\sC$
of $_c\Gamma_A$ can be obtained from indecomposable modules in stable tubes of concealed canonical algebras by a finite sequence
of admissible operations of types (ad~1)-(ad~5) and their duals (ad~1$^*$)-(ad~5$^*$) (see \cite[Section 3]{MS2} for details).
\end{remark}
We would like to stress that the cycle-finiteness assumption imposed on the infinite component $\sC$ of $_c\Gamma_A$ is essential for the validity
of the above theorem. Namely, it has been proved in \cite{Sk11}, \cite{Sk12} that, for an arbitrary finite dimensional algebra $B$ over a field
$K$, a module $M$ in $\mo B$, and a positive integer $r$, there exists a finite dimensional algebra $A$ over $K$ such that $B$ is a quotient
algebra of $A$, $\Gamma_A$ admits a faithful generalized standard stable tube $\sT$ of rank $r$, $\sT$ is not cycle-finite, and $M$
is a subfactor of all but finitely many indecomposable modules in $\sT$. This shows that in general the problem of describing the support
algebras of infinite cyclic components (even stable tubes) of Auslander-Reiten quivers is difficult.

\bigskip

In order to present the second case (when a cycle-finite cyclic component $\sC$ is finite), we need the class of generalized double tilted algebras
introduced by I. Reiten and A. Skowro\'nski in \cite{RS2} (see also \cite{AC}, \cite{CL} and \cite{RS1}).
A \textit{generalized double tilted algebra} is an algebra $B$ for which $\Gamma_B$ admits a separating almost acyclic component $\sC$.

The following consequence of \cite[Section 3]{RS2} describes the structure of the Auslander-Reiten quivers of generalized double tilted algebras.

\begin{theorem} \label{thm-45}
Let $B$ be a generalized double tilted algebra. Then the Auslander-Reiten quiver $\Gamma_B$ has a disjoint union decomposition
$\Gamma_B={\sP}^B \cup {\sC}^B \cup {\sQ}^B$, where
\begin{enumerate}
\item[\rm(1)] ${\sC}^B$ is an almost acyclic component separating ${\sP}^B$ from ${\sQ}^B$;
\item[\rm(2)] For each $i\in\{1,\ldots, m\}$, there exists hereditary algebra $H_i^{(l)}$ and tilting module $T_i^{(l)}\in\mo H_i^{(l)}$
such that the tilted algebra $B_i^{(l)}$ $=$ $\End_{H_i^{(l)}}(T_i^{(l)})$ is a quotient algebra of $B$ and ${\sP}^B$ is the
disjoint union of all components of $\Gamma_{B_i^{(l)}}$ contained entirely in the torsion-free part
${\mathscr Y}(T_i^{(l)})$ of $\mo B_i^{(l)}$ determined by $T_i^{(l)}$;
\item[\rm(3)] For each $j\in\{1,\ldots, n\}$, there exists hereditary algebra $H_j^{(r)}$ and tilting module $T_j^{(r)}\in\mo H_j^{(r)}$
such that the tilted algebra $B_j^{(r)}$ $=$ $\End_{H_j^{(r)}}(T_j^{(r)})$ is a quotient algebra of $B$ and ${\sQ}^B$ is the
disjoint union of all components of $\Gamma_{B_j^{(r)}}$ contained entirely in the torsion part
${\mathscr X}(T_j^{(r)})$ of $\mo B_j^{(r)}$ determined by $T_j^{(r)}$;
\item[\rm(4)] Every indecomposable module in ${\sC}^B$ not lying in the core $c({\sC}^B)$ of ${\sC}^B$ is an indecomposable module
over one of the tilted algebras $B_1^{(l)}, \ldots, B_m^{(l)}$, $B_1^{(r)}, \ldots,$ $B_n^{(r)}$;
\item[\rm(5)] Every nondirecting indecomposable module in ${\sC}^B$ is cycle-finite and lies in $c({\sC}^B)$;
\item[\rm(6)] $\pd_BX\leq 1$ for all indecomposable modules $X$ in ${\sP}^B$;
\item[\rm(7)] $\id_BY\leq 1$ for all indecomposable modules $Y$ in ${\sQ}^B$;
\item[\rm(8)] For all but finitely many indecomposable modules $M$ in ${\sC}^B$, we have $\pd_BM\leq 1$ or $\id_BM\leq 1$.
\end{enumerate}
\end{theorem}

Then ${\sC}^B$ is called a \textit{connecting component} of $\Gamma_B$, $B^{(l)}=B_1^{(l)}\times\ldots\times B_m^{(l)}$ is called the \textit{left tilted algebra} of $B$,
and $B^{(r)}=B_1^{(r)}\times\ldots\times B_n^{(r)}$ is called the \textit{right tilted algebra} of $B$.
Further, the almost acyclic component ${\sC}^B$ of $B$ admits a multisection.
Recall that, following \cite[Section 2]{RS2},
a full connected subquiver $\Delta$ of $\sC$ is called a \textit{multisection} if the following conditions are satisfied:
\begin{enumerate}
\item[\rm(1)] $\Delta$ is almost acyclic.
\item[\rm(2)] $\Delta$ is convex in $\sC$.
\item[\rm(3)] For each $\tau_B$-orbit $\mathcal O$ in $\sC$, we have $1 \leq | \Delta \cap {\mathcal O} | < \infty$.
\item[\rm(4)] $| \Delta \cap {\mathcal O} | = 1$ for all but finitely many $\tau_B$-orbits $\mathcal O$ in $\sC$.
\item[\rm(5)] No proper full convex subquiver of $\Delta$ satisfies \rm{(1)}--\rm{(4)}.
\end{enumerate}

Moreover, for a multisection $\Delta$ of a component $\sC$, the following full subquivers of $\sC$ were defined in \cite{RS2}:
$$\! \Delta^{\prime}_{l} = \{ X\! \in \!\Delta; \text{there is a nonsectional path in $\sC$ from $X$ to a projective module $P$}\},$$
$$\! \Delta^{\prime}_{r} = \{ X\! \in \!\Delta; \text{there is a nonsectional path in $\sC$ from an injective module $I$ to $X$}\},$$
$$
  \Delta^{\prime\prime}_{l} =
  \{ X \in \Delta^{\prime}_{l};
     \tau_A^{-1} X \notin \Delta^{\prime}_{l}
  \} , \qquad
  \Delta^{\prime\prime}_{r} =
  \{ X \in \Delta^{\prime}_{r};
     \tau_A X \notin \Delta^{\prime}_{r}
  \} ,
$$
$$
  \Delta_l =
  (\Delta \setminus \Delta^{\prime}_r)
       \cup \tau_A \Delta^{\prime\prime}_r
    , \quad
  \Delta_c = \Delta^{\prime}_l \cap \Delta^{\prime}_r
    , \quad
  \Delta_r = (\Delta \setminus \Delta^{\prime}_l)
        \cup \tau_A^{-1} \Delta^{\prime\prime}_l
  .
$$

\noindent Then $\Delta_l$ is called the {\it left part} of $\Delta$, $\Delta_r$ the {\it right part} of $\Delta$, and $\Delta_c$ the {\it core} of $\Delta$.
The following basic properties of $\Delta$ have been established in \cite[Proposition 2.4]{RS2}:
\begin{itemize}
\item Every cycle of $\sC$ lies in $\Delta_c$.
\item $\Delta_c$ is finite.
\item Every indecomposable module $X$ in $\sC$ is in $\Delta_c$, or a predecessor of $\Delta_l$ or a successor of $\Delta_r$ in $\sC$.
\end{itemize}
\begin{remark}
The class of algebras of finite representation type coincides with the class of generalized double tilted algebras $B$ with $\Gamma_B$
being the connecting component ${\sC}^B$ (equivalently, with the tilted algebras $B^{(l)}$ and $B^{(r)}$ being of finite representation type (possibly empty)).
\end{remark}

The following theorem from \cite[Theorem 1.2]{MPS3} describes the support algebras of finite cycle-finite cyclic components of the Auslander-Reiten quivers
of artin algebras.

\begin{theorem} \label{thm-1-2}
Let $A$ be an algebra and $\sC$ be a finite cycle-finite component of $_c{\Gamma_A}$. Then the following statements hold.
\begin{enumerate}
\item[\rm(1)] $\Supp(\sC)$ is a generalized double tilted algebra.
\item[\rm(2)] $\sC$ is the core $c({\sC}^{B(\sC)})$ of the unique almost acyclic connecting component ${\sC}^{B(\sC)}$ of $\Gamma_{B(\sC)}$.
\item[\rm(3)] $\Supp(\sC)=B(\sC)$.
\end{enumerate}
\end{theorem}
\begin{remark}
It follows from \cite[Corollary 2.7]{MPS3} that every finite cyclic component $\sC$ of an Auslander-Reiten quiver $\Gamma_A$ contains both a projective module and
an injective module, and hence $\Gamma_A$ admits at most finitely many finite cyclic components.
\end{remark}

We may summarize the description of the support algebras of cycle-finite cyclic components of the Auslander-Reiten quivers of artin algebras as follows:
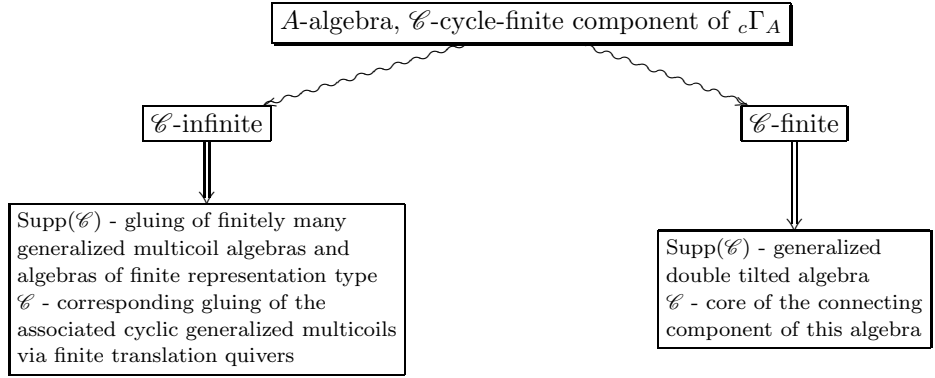
\[
\xymatrix@C=-50pt{
&*+[F-,]{\mbox{$A$-algebra, $\sC$-cycle-finite component of $_c\Gamma_A$}}\ar@{~>}[ld]\ar@{~>}[rd] \\
*+[F]{\mbox{$\sC$-infinite}}\ar@{=>}[d]&&*+[F]{\mbox{$\sC$-finite}}\ar@{=>}[d]\\
*+[F]{\shortstack[l]{\footnotesize$\Supp(\sC)$ - gluing of finitely many \\
\footnotesize generalized multicoil algebras and \\
\footnotesize algebras of finite representation type \\
\footnotesize $\sC$ - corresponding gluing of the \\
\footnotesize associated cyclic generalized multicoils \\
\footnotesize via finite translation quivers}}&&
*+[F]{\shortstack[l]{\footnotesize$\Supp(\sC)$ - generalized \\
\footnotesize double tilted algebra \\
\footnotesize $\sC$ - core of the connecting \\
\footnotesize component of this algebra}} \\
}
\]

Recall that an idempotent $e$ of an algebra $A$ is called \textit{convex} provided $e$ is a sum of pairwise orthogonal primitive idempotents of $A$
corresponding to the vertices of a convex valued subquiver of the quiver $Q_A$ of $A$. The following direct
consequence of Theorems \ref{thm-1-1}, \ref{thm-1-2} and \cite[Propositions 2.3, 2.4]{MPS3} provides a handy description of the faithful
algebra of a cycle-finite component of $_c\Gamma_A$.

\begin{corollary} \label{cor-50}
Let $A$ be an algebra and $\sC$ be a cycle-finite component of $_c\Gamma_A$. Then there exists
a convex idempotent $e_{\sC}$ of $A$ such that $\Supp(\sC)$ is isomorphic to the algebra $e_{\sC}Ae_{\sC}$.
\end{corollary}

\section{Cycle-finite indecomposable modules} \label{sect5}

The aim of this section is to present some results describing homological properties of indecomposable cycle-finite modules over artin algebras.

Let $A$ be an algebra and $M$ a~module in $\mo A$. We denote by $|M|$ the length of $M$ over the commutative artin ring $K$.
The following theorem is a consequence of Theorems \ref{thm-1-1}, \ref{thm-1-2}, the results established in \cite[Theorem 1.3]{MS3}
and the properties of directing modules described in \cite[2.4(8)]{Ri1}.
\begin{theorem} \label{thm-51}
Let $A$ be an algebra. Then, for all but finitely many isomorphism classes of cycle-finite modules $M$ in $\ind A$,
the following statements hold.
\begin{enumerate}
\item[\rm(1)] $|\Ext_A^1(M,M)|\leq |\End_A(M)|$ and $\Ext_A^r(M,M)=0$ for $r\geq 2$.
\item[\rm(2)] $|\Ext_A^1(M,M)| = |\End_A(M)|$ if and only if there is a quotient concealed canonical algebra $C$ of $A$ and a stable tube $\sT$ of $\Gamma_C$
such that $M$ is an indecomposable $C$-module in $\sT$ of quasi-length divisible by the rank of $\sT$.
\end{enumerate}
\end{theorem}

In particular, the above theorem shows that, for all but finitely many isomorphism classes of cycle-finite modules $M$ in $\ind A$,
the Euler characteristic
\[ \chi_A(M)=\sum_{i=0}^{\infty}(-1)^i|\Ext_A^i(M,M)| \]
of $M$ is well defined and nonnegative.

Let $A$ be an algebra and $K_0(A)$ the Grothendieck group of $A$. For a module $M$ in $\mo A$, we denote by $[M]$ the image of $M$ in $K_0(A)$.
Then $K_0(A)$ is a free abelian group with a ${\Bbb Z}$-basis given by $[S_1], \ldots, [S_n]$ for a complete family $S_1, \ldots, S_n$ of pairwise
nonisomorphic simple modules in $\mo A$. Thus, for modules $M$ and $N$ in $\mo A$, we have $[M]=[N]$ if and only if the modules $M$ and $N$ have
the same composition factors including the multiplicities. In particular, it would be interesting to find sufficient conditions for a module $M$
in $\ind A$ to be uniquely determined (up to isomorphism) by its composition factors (see \cite{RSS} for a general result in this direction).

The next theorem provides information on the composition factors of cycle-finite modules, and is a direct consequence
of Theorems \ref{thm-1-1}, \ref{thm-1-2}, \ref{thm-51} and the results established in \cite[Theorems A and B]{Ma}.
\begin{theorem} \label{thm-55}
Let $A$ be an algebra. The following statements hold.
\begin{enumerate}
\item[\rm(1)] There is a positive integer $m$ such that, for any cycle-finite module $M$ in $\ind A$ with $|\End_A(M)| \neq |\Ext_A^1(M,M)|$, the number
of isomorphism classes of modules $X$ in $\ind A$ with $[X]=[M]$ is bounded by $m$.
\item[\rm(2)] For all but finitely many isomorphism classes of cycle-finite modules $M$ in $\ind A$ with $|\End_A(M)| = |\Ext_A^1(M,M)|$, there are infinitely
many pairwise nonisomorphic modules $X$ in $\ind A$ with $[X]=[M]$.
\end{enumerate}
\end{theorem}

Following M. Auslander and I. Reiten \cite{AR}, one associates with each nonprojective module $X$ in $\ind A$ the number $\alpha(X)$
of indecomposable direct summands in the middle term
\[ 0\to \tau_AX\to Y\to X\to 0 \]
of an almost split sequence with the right term $X$. It has been proved by R. Bautista and S. Brenner \cite{BaBr} that, if $A$ is an algebra of finite
representation type and $X$ is a nonprojective module in $\ind A$, then $\alpha(X)\leq 4$, and if $\alpha(X)=4$ then $Y$ admits a projective-injective
indecomposable direct summand $P$, and hence $X=P/\soc(P)$.
In \cite{Li4} S. Liu proved that the same is true for any indecomposable nonprojective module $X$ lying on an oriented cycle of the Auslander-Reiten
quiver $\Gamma_A$ of any algebra $A$, and consequently for any nonprojective and nondirecting cycle-finite module in $\ind A$.

The following theorem is a direct consequence of Theorems \ref{thm-1-1}, \ref{thm-1-2}, and \cite[Corollary B]{MS1},
and provides more information on almost split sequences of nondirecting cycle-finite modules.

\begin{theorem} \label{thm-57}
Let $A$ be an algebra. Then, for all but finitely many isomorphism classes of nonprojective and nondirecting cycle-finite modules $M$
in $\ind A$, we have $\alpha(M)\leq 2$.
\end{theorem}

\section{Cycle-finite artin algebras} \label{sect6}

Following I. Assem and A. Skowro\'nski \cite{AS1}, \cite{AS2}, an algebra $A$ is said to be \textit{cycle-finite} if all cycles in $\ind A$ are finite.
A generalized multicoil algebra $B$ is called \textit{tame} if the quasitilted algebras $B^{(l)}$ and $B^{(r)}$ (see Section \ref{sect4} for definitions)
are products of tilted algebras of Euclidean types or tubular algebras.
Moreover, a generalized double tilted algebra $C$ is called \textit{tame} if the tilted algebras $C^{(l)}$ and $C^{(r)}$ (see Section \ref{sect4} for definitions)
are generically tame (see Section \ref{sect7} for definition) in the sense of W. Crawley-Boevey \cite{CB1}, \cite{CB2}.
We note that every tame generalized multicoil algebra and every tame generalized double tilted algebra is a cycle-finite algebra.

The following theorem describe the structure of the category $\ind A$ of an arbitrary cycle-finite algebra $A$, and
is a direct consequence of Theorems \ref{thm-1-1} and \ref{thm-1-2} (see also \cite[Theorems 7.1, 7.2 and 7.3]{MPS1}).

\begin{theorem} \label{thm7}
Let $A$ be a cycle-finite algebra. Then there exist tame generalized multicoil algebras $B_1, \ldots, B_p$ and tame
generalized double tilted algebras $B_{p+1}, \ldots, B_q$ which are quotient algebras of $A$ and the following statements hold.
\begin{enumerate}
\item[\rm(1)] $\ind A = \bigcup_{i=1}^q\ind B_i$.
\item[\rm(2)] All but finitely many isomorphism classes of modules in $\ind A$ belong to $\bigcup_{i=1}^p\ind B_i$.
\item[\rm(3)] All but finitely many isomorphism classes of nondirecting modules in $\ind A$ belong to generalized multicoils of $\Gamma_{B_1}, \ldots, \Gamma_{B_p}$.
\end{enumerate}
\end{theorem}

The next theorem extends the homological characterization of strongly simply connected algebras of polynomial growth established in \cite{PS1} to arbitrary
cycle-finite algebras, and is a direct consequence of Theorem \ref{thm-51} and the properties of directing modules described in \cite[2.4(8)]{Ri1}.

\begin{theorem} \label{thm-62}
Let $A$ be a cycle-finite algebra. Then, for all but finitely many isomorphism classes of modules $M$ in $\ind A$,
we have $|\Ext_A^1(M,M)|\leq |\End_A(M)|$ and $\Ext_A^r(M,M)$ $=0$ for $r\geq 2$.
\end{theorem}

In connection to Theorem \ref{thm-57} we present the following theorem proved by J. A. de la Pe\~na and the authors in \cite[Main Theorem]{MPS2}.

\begin{theorem} \label{thm-69}
Let $A$ be a cycle-finite algebra and $M$ be a nonprojective module in $\ind A$, and
\[ 0\to \tau_AM\to N\to M\to 0\]
be the associated almost split sequence in $\mo A$. The following statements hold.
\begin{enumerate}
\item[\rm(1)] $\alpha(M)\leq 5$.
\item[\rm(2)] If $\alpha(M)=5$ then $N$ admits an indecomposable projective-injective direct summand $P$, and hence $M\simeq P/\soc(P)$.
\end{enumerate}
\end{theorem}

The following example of a~cycle-finite algebra $A$ from \cite[Example 5.17]{Ma11}, illustrates the above theorem.

\begin{example} \label{ex-5-cf}
Let $K$ be a field and $A=KQ/I$ the bound quiver algebra over $K$, where $Q$ is the quiver of the form
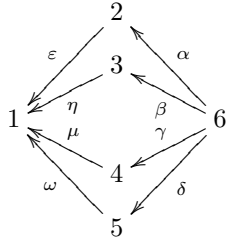
\[ \xymatrix@R=8pt@C=28pt{
&2\ar[ldd]_{\varepsilon}\\
&3\ar[ld]^{\eta}\\
1&&6\ar[luu]_{\alpha}\ar[lu]^{\beta}\ar[ld]_{\gamma}\ar[ldd]^{\delta}\\
&4\ar[lu]_{\mu}\\
&5\ar[luu]^{\omega}\\
} \]
and $I$ the ideal of the path algebra $KQ$ of $Q$ over $K$ generated by the paths
$\beta\eta - \alpha\varepsilon$, $\gamma\mu - \alpha\varepsilon$, $\delta\omega - \alpha\varepsilon$.
Denote by $B$ the hereditary algebra given by the full subquiver of $Q$ given by the vertices $1, 2, 3, 4, 5$ and by $C$
the hereditary algebra given by the full subquiver of $Q$ given by the vertices $2, 3, 4, 5, 6$.
Note that $P_6 = I_1$ is projective-injective $A$-module. Therefore, applying \cite[Proposition IV.3.11]{ASS}, we conclude that there is
in $\mo A$ an almost split sequence of the form
\[ 0 \longrightarrow \rad P_6 \longrightarrow S_2\oplus S_3\oplus S_4\oplus S_5\oplus P_6 \longrightarrow P_6/S_1 \longrightarrow 0, \]
where $S_2\oplus S_3\oplus S_4\oplus S_5 \cong \rad P_6/S_1$.
Moreover, $\rad P_6$ is the indecomposable injective $B$-module $I_1^B$, whereas $P_6/S_1$ is the indecomposable projective $C$-module
$P_6^C$. The component of $\Gamma_A$ containing $P_6 = I_1$ is the following gluing of the preinjective component of $\Gamma_B$ with
the postprojective component of $\Gamma_C$ (see details in \cite[Example VIII.5.7(e)]{ASS})
\[ \xymatrix@R=10pt@C=14pt{
\ldots{\phantom{P}}\ar[rdd]&&\tau_BS_2\ar[rdd]&&S_2\ar[rdd]&&\tau_C^{-1}S_2\ar[rdd]&&{\phantom{P}}\ldots \\
\ldots{\phantom{P}}\ar[rd]&&\tau_BS_3\ar[rd]&&S_3\ar[rd]&&\tau_C^{-1}S_3\ar[rd]&&{\phantom{P}}\ldots \\
&\tau_BI_1^B\ar[ruu]\ar[ru]\ar[rd]\ar[rdd]&&I_1^B\ar[r]\ar[ruu]\ar[ru]\ar[rd]\ar[rdd]&P_6\ar[r]&P_6^C\ar[ruu]\ar[ru]\ar[rd]\ar[rdd]&&
\tau_C^{-1}P_6^C\ar[ruu]\ar[ru]\ar[rd]\ar[rdd] \\
\cdots{\phantom{P}}\ar[ru]&&\tau_BS_4\ar[ru]&&S_4\ar[ru]&&\tau_C^{-1}S_4\ar[ru]&&{\phantom{P}}\cdots \\
\cdots{\phantom{P}}\ar[ruu]&&\tau_BS_5\ar[ruu]&&S_5\ar[ruu]&&\tau_C^{-1}S_5\ar[ruu]&&{\phantom{P}}\cdots \\
} \]
where $\tau_B$ and $\tau_C$ denote the Auslander-Reiten translations in $\mo B$ and $\mo C$, respectively.
\end{example}

Let $A$ be an algebra. Recall that, following C. M. Ringel \cite{Ri1}, a module $M$ in $\ind A$ which does not lie on a cycle in
$\ind A$ is called \textit{directing}.
We note that if all modules in $\ind A$ are directing, then $A$ is of finite representation type \cite{Ri1}
(see also \cite{HL} for the corresponding result over arbitrary artin algebra).
If $A$ is a cycle-finite algebra and $M$ is a module in $\ind A$, then $M$
is a directing if and only if $M$ is an acyclic vertex of $\Gamma_A$. The following result from \cite{BiSk}
provides solution of the open problem concerning infinity of directing modules for cycle-finite algebras of infinite representation type.

\begin{theorem} \label{thm-65}
Let $A$ be a cycle-finite algebra of infinite representation type. Then $\ind A$ contains infinitely many directing modules.
\end{theorem}

In connection to the above theorem the following question occurred naturally:
does every cycle-finite algebra of infinite representation type contains at least one directing projective module
or a directing injective module?

In general, the answer is negative. Namely, in the joint work with J. A. de la Pe\~na \cite{MPS4} we
constructed a family of cycle-finite algebras of infinite representation type
with all indecomposable projective modules and indecomposable injective modules nondirecting.
Moreover, it has been shown (see \cite[Theorem]{MPS4}) that there are such algebras with an arbitrary large number of almost acyclic
Auslander-Reiten components having finite cyclic multisections.

The following example from \cite[Section 9]{MPS4}, illustrates the situation.
We refer to \cite{MPS4} for the general construction of algebras of this type and more examples.
\begin{example} \label{ex-cm2017}
Let $K$ be an algebraically closed field and $A=K\Sigma/J$ the bound quiver algebra over $K$, where $\Sigma$ is the quiver of the form
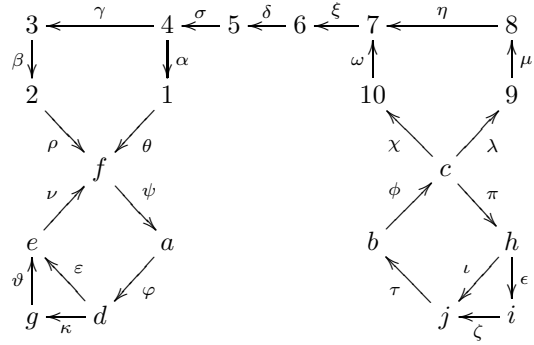
\[
\xymatrix@R=14pt@C=14pt{
3\ar[d]_{\beta}&&4\ar[ll]_{\gamma}\ar[d]^{\alpha}&5\ar[l]_{\sigma}&6\ar[l]_{\delta}&7\ar[l]_{\xi}&&8\ar[ll]_{\eta}\cr
2\ar[rd]_{\rho}&&1\ar[ld]^{\theta}&&&10\ar[u]^{\omega}&&9\ar[u]_{\mu}\cr
&f\ar[rd]^{\psi}&&&&&c\ar[lu]^{\chi}\ar[ru]_{\lambda}\ar[rd]^{\pi}\cr
e\ar[ru]^{\nu}&&a\ar[ld]^{\varphi}&&&b\ar[ru]^{\phi}&&h\ar[ld]_{\iota}\ar[d]^{\epsilon}\cr
g\ar[u]^{\vartheta}&d\ar[lu]_{\varepsilon}\ar[l]^{\kappa}&&&&&j\ar[lu]^{\tau}&i\ar[l]^{\zeta}\cr
}
\]
and $J$ is the ideal in $K\Sigma$ generated by the elements $\omega\xi\delta\sigma\alpha$, $\alpha\theta-\gamma\beta\rho$, $\chi\omega-\lambda\mu\eta$, $\rho\psi$, $\theta\psi$,
$\nu\psi$, $\varphi\varepsilon$, $\varphi\kappa$, $\kappa\vartheta$, $\vartheta\nu$, $\phi\chi$, $\phi\lambda$, $\phi\pi$, $\pi\epsilon$, $\epsilon\zeta$, $\zeta\tau$, $\iota\tau$.
Denote by $B=KQ/I$ the bound quiver algebra given by the subquiver $Q$ of $\Sigma$ given by the vertices $1, 2, \ldots, 10$ and the ideal $I$ in $KQ$
generated by the element $\omega\xi\delta\sigma\alpha$. Then $B$ is a~tubular algebra of type $(2,3,6)$ in the sense of Ringel \cite{Ri1} and following \cite[Section 7]{MPS4}
the algebra $B$ is called \textit{exceptional tubular algebra of type $(2,3,6)$}.
Moreover, we denote by $H_0 = KQ^{(0)}$ the hereditary algebra of Euclidean type $\widetilde{{\Bbb E}}_8$ given by the full subquiver $Q^{(0)}$ of $\Sigma$
with the vertices $1, 2, 3, 4, 5, 6, 7, 8, 9$, and by $H_1 = KQ^{(1)}$ the hereditary algebra of Euclidean type $\widetilde{{\Bbb E}}_8$ given by the full subquiver
$Q^{(1)}$ of $\Sigma$ with the vertices $2, 3, 4, 5, 6, 7, 8, 9, 10$.

Then it follows from \cite[Theorems 5.1, 5.2 and 7.1]{MPS4} that the Auslander-Reiten quiver $\Gamma_A$ of $A$ has a~decomposition
\[ \Gamma_A = {\sC}_0 \cup {\sT}^B_0 \cup \bigg(\bigcup_{q\in{\Bbb Q}_1^0}{\sT}^B_q\bigg) \cup {\sT}^B_{1} \cup {\sC}_1, \]
where
\begin{itemize}
\item ${\sC}_0$ is an almost acyclic component of the form $\Delta^{(0)}\cup{\sP}^B$ with a~cyclic multisection $\Delta^{(0)}$ such that the left part $\Delta^{(0)}_l$
of $\Delta^{(0)}$ is empty, the right part $\Delta^{(0)}_r$ of $\Delta^{(0)}$ is given by the indecomposable projective modules of the postprojective component ${\sP}^B$ of
$\Gamma_B$, being the postprojective component ${\sP}^{H_0}$,
\item ${\sT}^B_0$ is the ${\Bbb P}_1(K)$-family $({\sT}^B_{0,\lambda})_{\lambda\in{\Bbb P}_1(K)}$ of ray tubes obtained from the ${\Bbb P}_1(K)$-family
${\sT}^{H_0}$ of stable tubes of tubular type $(2,3,5)$ by inserting in the unique stable tube of rank $5$, say ${\sT}^{H_0}_{0,0}$, one ray as follows
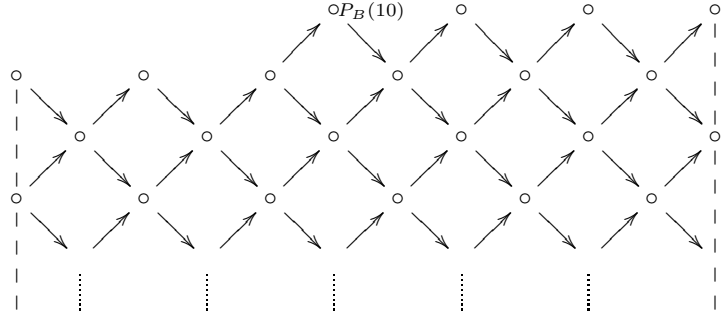
\[
\xymatrix@C=13pt@R=13pt{
&&&&&\circ\ar[rd]&\hspace{-7mm}\scr{P_B(10)}&\circ\ar[rd]&&\circ\ar[rd]&&\circ\ar@{--}[dd]\cr
\circ\ar@{--}[dd]\ar[rd]&&\circ\ar[rd]&&\circ\ar[ru]\ar[rd]&&\circ\ar[ru]\ar[rd]&&\circ\ar[ru]\ar[rd]&&\circ\ar[ru]\ar[rd]\cr
&\circ\ar[ru]\ar[rd]&&\circ\ar[ru]\ar[rd]&&\circ\ar[ru]\ar[rd]&&\circ\ar[ru]\ar[rd]&&\circ\ar[ru]\ar[rd]&&\circ\ar@{--}[ddd]\cr
\circ\ar[ru]\ar[rd]\ar@{--}[dd]&&\circ\ar[ru]\ar[rd]&&\circ\ar[ru]\ar[rd]&&\circ\ar[ru]\ar[rd]&&\circ\ar[ru]\ar[rd]&&\circ\ar[ru]\ar[rd]\cr
{\phantom{\circ}}&{\phantom{\circ}}\ar[ru]\ar@{.}[d]&&{\phantom{\circ}}\ar[ru]\ar@{.}[d]&&{\phantom{\circ}}\ar[ru]\ar@{.}[d]&&{\phantom{\circ}}\ar[ru]\ar@{.}[d]&&{\phantom{\circ}}\ar[ru]\ar@{.}[d]&&{\phantom{\circ}}\cr
{\phantom{\circ}}&{\phantom{\circ}}&&{\phantom{\circ}}&&{\phantom{\circ}}&&{\phantom{\circ}}&&{\phantom{\circ}}&&{\phantom{\circ}}\cr
}
\]
where the vertices along to the dashed vertical lines have to be identified,
\item For each $q\in{\Bbb Q}_1^0$, ${\sT}^B_q$ is a ${\Bbb P}_1(K)$-family $({\sT}^B_{q,\lambda})_{\lambda\in{\Bbb P}_1(K)}$ of stable tubes of tubular
type $(2,3,6)$, where ${\Bbb Q}_1^0 = {\Bbb Q}\cap (0,1)$,
\item ${\sT}^B_{1}$ is the ${\Bbb P}_1(K)$-family $({\sT}^B_{1,\lambda})_{\lambda\in{\Bbb P}_1(K)}$ of coray tubes obtained from the
${\Bbb P}_1(K)$-family ${\sT}^{H_{1}}$ of stable tubes of tubular type $(2,3,5)$ by inserting in the unique stable tube of rank $5$,
say ${\sT}^{H_{1}}_{1,0}$, one coray as follows
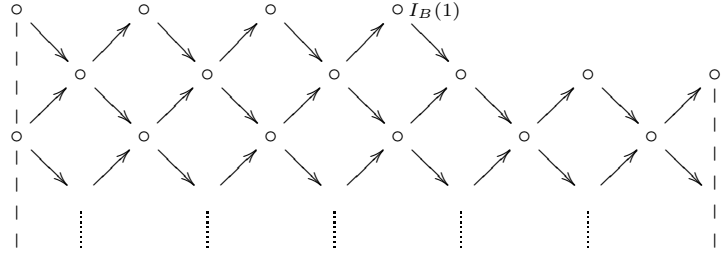
\[
\xymatrix@C=13pt@R=13pt{
\circ\ar@{--}[dd]\ar[rd]&&\circ\ar[rd]&&\circ\ar[rd]&&\circ\ar[rd]&\hspace{-7mm}\scr{I_B(1)}\cr
&\circ\ar[ru]\ar[rd]&&\circ\ar[ru]\ar[rd]&&\circ\ar[ru]\ar[rd]&&\circ\ar[rd]&&\circ\ar[rd]&&\circ\ar@{--}[ddd]\cr
\circ\ar[ru]\ar[rd]\ar@{--}[dd]&&\circ\ar[ru]\ar[rd]&&\circ\ar[ru]\ar[rd]&&\circ\ar[ru]\ar[rd]&&\circ\ar[ru]\ar[rd]&&\circ\ar[ru]\ar[rd]\cr
{\phantom{\circ}}&{\phantom{\circ}}\ar[ru]\ar@{.}[d]&&{\phantom{\circ}}\ar[ru]\ar@{.}[d]&&{\phantom{\circ}}\ar[ru]\ar@{.}[d]&&{\phantom{\circ}}\ar[ru]\ar@{.}[d]&&{\phantom{\circ}}\ar[ru]\ar@{.}[d]&&{\phantom{\circ}}\cr
{\phantom{\circ}}&{\phantom{\circ}}&&{\phantom{\circ}}&&{\phantom{\circ}}&&{\phantom{\circ}}&&{\phantom{\circ}}&&{\phantom{\circ}}\cr
}
\]
where the vertices along to the dashed vertical lines have to be identified,
\item ${\sC}_1$ is an almost acyclic component of the form ${\sQ}^B\cup\Delta^{(1)}$ with a~cyclic multisection $\Delta^{(1)}$ such that the right part $\Delta^{(1)}_r$
of $\Delta^{(1)}$ is empty, the left part $\Delta^{(1)}_l$ of $\Delta^{(1)}$ is given by the indecomposable injective modules of the preinjective component ${\sQ}^B$ of
$\Gamma_B$, being the preinjective component ${\sQ}^{H_1}$.
\end{itemize}
Therefore, all indecomposable projective modules and indecomposable injective modules in $\mo A$ are nondirecting.
\end{example}

Let $A$ be an algebra.
We denote by $\cF(A)$ the category of all finitely presented contravariant functors from $\mo A$ to the category $\Ab$
of abelian groups. The category $\cF(A)$ was intensively studied over the last 40 years, and is considered to be one
of the important topics of the modern representation theory of algebras.
It is a hard problem to describe the category $\cF(A)$ even if the category $\mo A$ is well understood.
A natural approach to study the structure of $\cF(A)$ is via the associated Krull-Gabriel filtration
\[
  0 = \cF(A)_{-1}
    \subseteq \cF(A)_0
    \subseteq \cF(A)_1
    \subseteq \dots
    \subseteq \cF(A)_{n-1}
    \subseteq \cF(A)_n
    \subseteq \cdots
\]
of $\cF(A)$ by Serre subcategories, where, for each $n \in \bN$, $\cF(A)_n$
is the subcategory of all functors $F$ in $\cF(A)$ which become of finite length in the quotient
category $\cF(A)/\cF(A)_{n-1}$ \cite{F}, \cite{Po}.
Following W. Geigle \cite{Ge1}, we define
$KG(A) = \min \{ n \in \bN \,|\, \cF(A)_n = \cF(A) \}$
if such a minimum exists, and set $KG(A) = \infty$ if it is not the case.
Then $KG(A)$ is called the \textit{Krull-Gabriel dimension} of $A$.
The interest in the Krull-Gabriel dimension $KG(A)$ is motivated by the fact that the above
filtration of $\cF(A)$ leads to a hierarchy of exact sequences in $\mo A$, where the almost
split sequences form the lowest level (see \cite{Ge1}).

The following characterization of cycle-finite algebras with finite Krull-Gabriel dimension,
has been established by the second named author in \cite[Theorem 1.2]{Sk-20}.

\begin{theorem} \label{thm-KG}
Let $A$ be a cycle-finite algebra of infinite representation type. The following statements are equivalent.
\begin{enumerate}
 \item[\rm(1)] $KG(A) < \infty$.
 \item[\rm(2)] $KG(A) = 2$.
 \item[\rm(3)] $\bigcap_{m \geq 1} (\rad_A^{\infty})^m = 0$.
 \item[\rm(4)] $\rad_A^{\infty}$ is nilpotent.
 \item[\rm(5)] All but finitely many components of $\Gamma_A$ are stable tubes of rank one.
 \item[\rm(6)] $A$ does not admit a tubular quotient algebra.
\end{enumerate}
\end{theorem}

We end this section with the related open problem (see \cite[Question 1]{MPS4}).
\begin{prob} \textit{
Let $A$ be a~cycle-finite algebra of infinite representation type and finite Krull-Gabriel dimension.
Is it true that $\ind A$ admits a~directing projective module or a~directing injective module?
}
\end{prob}

\section{Artin algebras with separating Auslander-Reiten components} \label{sect7}

In this section we discuss the structure of artin algebra $A$ having separating families of components in $\Gamma_A$.

Let $A$ be an algebra. A component $\sC$ of $\Gamma_A$ is called \textit{cycle-finite} if all modules in $\sC$ are cycle-finite.

We have the following proposition.

\begin{proposition} \label{pro-71}
Let $A$ be an algebra with a separating family ${\sC}^A$ of components in $\Gamma_A$, and $\Gamma_A$=${\sP}^A \cup {\sC}^A \cup {\sQ}^A$
the associated decomposition of $\Gamma_A$. Then ${\sC}^A$ is a family of cycle-finite components.
\end{proposition}

The next theorem from \cite[Theorem 1.5]{MS9} provides solution of the problem concerning the structure of artin algebras admitting a~separating family
of Auslander-Reiten components, initiated by Ringel \cite{Ri0}, \cite{Ri1}, \cite{Ri-86}.
\begin{theorem}\label{thm-75}
Let $A$ be an algebra with a separating family ${\sC}^A$ of components in $\Gamma_A$, and $\Gamma_A$=${\sP}^A \cup {\sC}^A \cup {\sQ}^A$
the associated decomposition of $\Gamma_A$. Then there exist quotient algebras $A^{(l)}$ and $A^{(r)}$ of $A$ such that the following statements hold.
\begin{enumerate}
\item[\rm(1)] $A^{(l)}=A^{(l)}_1 \times\cdots\times A^{(l)}_m \times A^{(l)}_{m+1} \times\cdots\times A^{(l)}_{m+p},$  where
\begin{enumerate}
\item[\rm(a)] For each $i\in\{1,\ldots,m\}$, $A^{(l)}_i$ is a~tilted algebra of the form $\End_{H^{(l)}_i}(T^{(l)}_i)$ for a~hereditary algebra $H^{(l)}_i$
and a~tilting module $T^{(l)}_i$ in $\mo H^{(l)}_i$ without indecomposable preinjective  direct summands.
\item[\rm(b)] For each $i\in\{m+1,\ldots,m+p\}$, $A^{(l)}_i$ is a~quasitilted algebra of canonical type with a~separating family of coray tubes in $\Gamma_{A^{(l)}_i}$.
\end{enumerate}
\item[\rm(2)] $A^{(r)}=A^{(r)}_1 \times\cdots\times A^{(r)}_n \times A^{(r)}_{n+1} \times\cdots\times A^{(r)}_{n+q},$ where
\begin{enumerate}
\item[\rm(a)] For each $j\in\{1,\ldots,n\}$, $A^{(r)}_j$ is a~tilted algebra of the form $\End_{H^{(r)}_j}(T^{(r)}_j)$ for a~hereditary algebra $H^{(r)}_j$
and a~tilting module $T^{(r)}_j$ in $\mo H^{(r)}_j$ without indecomposable postprojective direct summands.
\item[\rm(b)] For each $j\in\{n+1,\ldots,n+q\}$, $A^{(r)}_j$ is a~quasitilted algebra of canonical type with a~separating family of ray tubes in $\Gamma_{A^{(r)}_j}$.
\end{enumerate}
\item[\rm(3)]  $\sP^A = \bigcup_{i=1}^{m+p} \sP^{A^{(l)}_i}$
and every component in $\sP^A$ is either a~postprojective component, a~ray tube, or obtained from a~component of the form ${\Bbb Z}{\Bbb A}_{\infty}$
by a~finite number (possibly zero) of ray insertions.
\item[\rm(4)]  $\sQ^A = \bigcup_{j=1}^{n+q} \sQ^{A^{(r)}_j}$
and every component in $\sQ^A$ is either a~preinjective component, a~coray tube, or obtained from a~component of the form ${\Bbb Z}{\Bbb A}_{\infty}$
by a~finite number (possibly zero) of coray insertions.
\end{enumerate}
\end{theorem}

In \cite{CB1}, \cite{CB2} Crawley-Boevey introduced the concept of a~generically tame algebra. An indecomposable right $A$-module $M$ over an algebra $A$ is called
a \textit{generic module} if $M$ is of infinite length over $A$ but of finite length over $\End_A(M)$, called the \textit{endolength} of $M$. Then an algebra $A$ is called
\textit{generically tame} if, for any positive integer $d$, there are only finitely many isomorphism classes of generic right $A$-modules of endolength $d$.
An algebra $A$ is called \textit{generically finite} if there are at most finitely many pairwise non-isomorphic generic right $A$-modules.
Further, $A$ is called \textit{generically of polynomial growth} if there is a~positive integer $m$ such that for any positive integer $d$ the number of isomorphism
classes of generic right $A$-modules of endolength $d$ is at most $d^m$.
We note that every algebra $A$ of finite representation type is generically trivial, that is, there is no generic right $A$-module.
We also stress that by a theorem of Crawley-Boevey \cite[Theorem 4.4]{CB1}, if $A$ is an algebra over an algebraically closed field $K$,
then $A$ is generically tame if and only if $A$ is tame in the sense of Drozd \cite{Dro} (see also \cite{CB}, \cite{SS2}).

Recall also that following \cite{Sk7} the \textit{component quiver} $\Sigma_A$ of an algebra $A$ has the components of $\Gamma_A$ as vertices and
there is an arrow $\sC\to\sD$ in $\Sigma_A$ if $\rad_A^{\infty}(X,Y)\neq 0$, for some modules $X$ in $\sC$ and $Y$ in $\sD$.
In particular, a~component $\sC$ of $\Gamma_A$ is generalized standard if and only if there is no loop at $\sC$ in $\Sigma_A$.

The final theorem is a consequence of \cite[Theorem 1.10]{MS9}, and characterizes the cycle-finite algebras with separating families of Auslander-Reiten components.
\begin{theorem}\label{thm-77}
Let $A$ be an algebra with a~separating family of components in $\Gamma_A$. The following statements are equivalent:
\begin{enumerate}
\item[\rm(1)] $A$ is cycle-finite.
\item[\rm(2)] $A$ is generically tame.
\item[\rm(3)] $A$ is generically of polynomial growth.
\item[\rm(4)] $A^{(l)}$ and $A^{(r)}$ are products of tilted algebras of Euclidean type or tubular algebras.
\item[\rm(5)] $\Gamma_A$ is almost periodic.
\item[\rm(6)] $\Sigma_A$ is acyclic.
\end{enumerate}
\end{theorem}

%
%
%

\bibliographystyle{amsalpha}

\end{document}